\documentclass[12pt]{amsart}
\usepackage{amsmath,amsthm,amssymb}
\usepackage[matrix,arrow,curve]{xy}
\makeatletter\@addtoreset{equation}{section} \makeatother

\newtheorem{theorem}[equation]{Theorem}
\newtheorem{proposition}[equation]{Proposition}
\newtheorem{lemma}[equation]{Lemma}
\newtheorem{corollary}[equation]{Corollary}

\theoremstyle{definition}
\newtheorem{example}[equation]{Example}

\theoremstyle{remark}
\newtheorem{remark}[equation]{Remark}

\def \O {\mathcal{O}}
\def \H {\mathcal{H}}
\def \X {\mathcal{X}}
\def \Q {\mathbb{Q}}
\def \P {\mathbb{P}}
\def \F {\mathbb{F}}
\def \Z {\mathbb{Z}}
\def \C {\mathbb{C}}

\def \ge {\geqslant}
\def \le {\leqslant}

\def \phix {\varphi_{|-K_X|}}

\def \rk {\mathrm{rk}}
\def \mult {\mathrm{mult}}
\def \Bs {\mathrm{Bs}}
\def \Bir {\mathrm{Bir}}
\def \Aut {\mathrm{Aut}}
\def \Pic {\mathrm{Pic}}
\def \PSO {\mathrm{PSO}}

\def \t {~---\ }

\newcommand{\tit}{Birational rigidity and $\Q$-factoriality\\
of a singular double quadric}
\author{C.\,Shramov}
\title{\tit}
\thanks{The work was partially supported by RFFI grants
No. 04-01-00613 and No. 05-01-00353}
\email{shramov@mccme.ru}

\begin{document}

\maketitle

\begin{abstract}
We prove birational rigidity and calculate the group of birational 
automorphisms of a nodal $\Q$-factorial double cover $X$ of a smooth 
three-dimensional quadric branched over a quartic section. We also prove
that $X$ is $\Q$-factorial provided that it has at most $11$ singularities;
moreover, we give an example of a non-$\Q$-factorial variety of this type 
with $12$ simple double singularities.
\end{abstract}

\section{Introduction}
\label{section:setup}

One of the popular problems of birational geometry is 
to find all Mori fibrations (see~\cite{Matsuki} for
definition) that are birationally equivalent to a given Mori fibration
$\X\to T$, and to compute the group 
of birational automorphisms $\Bir(\X)$ of a variety $\X$. 
The cases when there are few structures of Mori fibrations on $\X$,
for example, when there is only one\t up to a natural equivalence 
(see below)\t structure of Mori fibration, are of special interest.

A Mori fibration $\X\to T$ is called \emph{birationally rigid} 
if for any birational map $\chi:\X\dasharrow\X'$ to another
Mori fibration $\X'\to T'$ there is a birational selfmap 
$\gamma:\X\dasharrow\X$ and a birational map $\sigma:T\dasharrow T'$
such that the following diagram commutes
$$
\xymatrix{
\X\ar@{-->}[r]^{\chi\circ\gamma}\ar[d] & \X'\ar[d]\\
T\ar@{-->}[r]^{\sigma} & T'
}
$$
Moreover, it is required that the map $\chi\circ\gamma$ restricts to 
an isomorphism on a general fiber.
Informally, this definition means that $\X$ cannot be transformed into 
a principally different Mori fibration.
A birationally rigid Fano variety $X$ with Picard number $\rho(X)=1$
such that its group of birational automorphisms
$\Bir(X)$ coincides with the subgroup $\Aut(X)$ of biregular ones is called
\emph{birationally superrigid}
(see~\cite{Corti} or~\cite{Pukhlikov-essentials}).

Fano varieties of low degree 
(a double cover of $\P^3$ branched over a sextic, a quartic and
a double cover of a quadric branched over a divisor of degree $4$\t
see~\cite{Iskovskikh-classification}) 
give examples of birationally rigid varieties
with relatively simple groups of birational selfmaps.
Birational superrigidity of a smooth quartic was proved 
in~\cite{IskovskikhManin}; 
a proof of birational superrigidity of a smooth double cover
of $\P^3$ branched over a sextic and birational rigidity of a
smooth double cover of a quadric branched over a divisor of degree
$4$ (together with the calculation of its group of birational automorphisms)
may be found in~\cite{Iskovskikh-rigid} and in~\cite{IskovskikhPukhlikov}.
The case of a quartic with one simple double point was considered 
in~\cite{Pukhlikov-quartic}
(where its birational rigidity was proved and a group of its birational 
automorphisms was calculated); the case of a $\Q$-factorial 
double cover of $\P^3$ branched over a sextic with arbitrary 
number of simple double points was studied in~\cite{CheltsovPark}
(this variety appeared to be birationally superrigid);
the case of a $\Q$-factorial quartic with arbitrary number of simple double
points was considered in~\cite{Mella}
(this variety was proved to be birationally rigid and the generators 
of its groups of birational automorphisms were found). 
In~\cite{Grinenko} birational rigidity of a double cover of a quadric branched
over a divisor of degree $4$ with one simple double point
was proved and the group of its birational selfmaps was calculated.
The main goal of this paper is to prove the following statement that
continues this series of results.

Let $Q\subset\P^4$ be a smooth quadric, $W\subset\P^4$\t 
a three-dimensional quartic, $S=W\cap Q$, $X$\t a double cover
of $Q$ branched over $S$; $X$ is a Fano variety of degree $\deg X=4$,
and the double cover structure is given by the map 
$\phix:X\to Q\subset\P^4$.

\begin{theorem}\label{maintheorem}
Assume that the variety $X$ is $\Q$-factorial
and has only simple double singularities. Then the following holds

1. $X$ is birationally rigid.

2. Let $B\subset X$\t be a line (i.\,e. a curve of anticanonical degree $1$)
such that $\phix(B)\not\subset S$.
Then there is a birational involution $\tau_B$ associated to $B$ 
(see Example~\ref{example:involution} for an explicit description); 
these involutions together with the subgroup $\Aut(X)$ generate the group
$\Bir(X)$.

3. Involutions $\tau_B$ are independent in the group $\Bir(X)$
(i.\,e. there are no relations on them, except for the trivial relations 
$\tau_B^2=1$).
Equivalently, there exists an exact sequence

$$
\xymatrix{
1\ar[r]&
F(\{\tau_B\})\ar[r]&
\Bir(X)\ar[r]&\Aut(X)\ar[r]&1
},
$$
where $F(\{\tau_B\})$ denotes a free group generated by all  
involutions $\tau_B$ associated to the lines $B$ such that  
$\phix(B)\not\subset S$.
\end{theorem}

\begin{remark}
The group $\Aut(X)$ of a general variety $X$ contains only one 
nontrivial element\t an involution $\delta$ of a double cover $\phix$.
\end{remark}

\begin{remark}
For $X$ to have only simple double singularities 
it is necessary and sufficient
that $S$ has only simple double singularities.
One may assume that $W$ is smooth along $S$ and is tangent to
$Q$ at the singular points of $S$.
\end{remark}

\begin{remark}
The assumption that $Q$ is smooth is important\t 
a double cover of a cone over a smooth quadric surface branched over 
a quartic section (not passing through a vertex) 
is a much more complicated variety. To start with, it is
not $\Q$-factorial. It has two small resolutions;
each of them covers a corresponding small resolution of
a quadric cone and is fibered into Del Pezzo surfaces of degree $2$
(these fibrations arise from two pencils of planes on a quadric cone).
On the other hand, all structures of Mori fibrations on this variety 
are these two (see~\cite{Grinenko-konus}).
\end{remark}

Theorem~\ref{maintheorem} is proved by the method of maximal 
singularities (see~\cite{Pukhlikov-essentials} or~\cite{Corti}). 
In section~\ref{section:preliminaries} we introduce 
the necessary notions (as well as the notations used throughout 
the paper), and after that we prove statements~1 and~2 
of Theorem~\ref{maintheorem} in sections~\ref{section:max-points} 
and~\ref{section:max-curves}.
The proof of statement~3 of Theorem~\ref{maintheorem} is contained in 
section~\ref{section:relations}.

Since Theorem~\ref{maintheorem} holds only for $\Q$-factorial varieties,
it seems natural to ask which conditions guarantee $\Q$-factoriality
of $X$. In section~\ref{section:Q-factoriality} 
we'll prove 
the following statement (note that an analogous statement was 
established for a sextic double solid in~\cite{CheltsovPark} and for a quartic 
in~\cite{Cheltsov-quartic}; see also~\cite{Shramov} for detailes about 
$\Q$-factorial quartics).

\begin{proposition}\label{proposition:Q-factoriality}
If the number of singular points of $X$ is at most~$11$, then $X$ is 
$\Q$-factorial. On the other hand, there exist non-$\Q$-factorial 
varieties of this type with $12$ simple double points.
\end{proposition}

All varieties are assumed to be defined over the field of complex numbers.
The degrees are always calculated with respect to the anticanonical 
linear system. We'll use the notations introduced in this section throughout 
the work.

\smallskip
The author is grateful to I.\,Cheltsov for numerous comments, 
and to S.\,Galkin and V.\,Przyjalkowski for 
useful conversations.

\section{Preliminaries on the method of maximal singularities}
\label{section:preliminaries}

In this section we briefly describe main constructions of the method
of maximal singularities and introduce the necessary 
notation (for a detailed treatment see~\cite{Pukhlikov-essentials} 
or~\cite{Corti}). All relevant definitions may be found, for example,
in~\cite{Matsuki}.

Let $V$ be a (three-dimensional) Fano variety with terminal 
singularities and Picard number $\rho(V)=1$, and let $V'\to S'$ be a 
Mori fibration (in a more general setting one may assume that $V$ is also 
a Mori fibration over an arbitrary base $S$, but we don't need it). 
Let $\chi:V\dasharrow V'$ be a rational map. There is an algorithm 
(known as Sarkisov program) that 
decomposes the map $\chi$ into a composition of elementary maps (links) 
of four types (see~\cite{Corti} or~\cite{Matsuki}). Choose a very ample divisor
$H'$ on $V'$ and let $\H=\chi_*^{-1}|H|'$ (note that $\H$ has no fixed 
components). Let $\mu$ be a (rational) number such that 
$\H\subset |-\mu K_V|$. The N\"oether--Fano inequality 
(see~\cite{Iskovskikh-rigid}, \cite{Corti}, \cite{Matsuki} 
or~\cite{Pukhlikov-essentials}) claims that if in this setting
$\chi$ is not an isomorphism then the pair 
$(V, \frac{1}{\mu}\H)$ is not canonical. Next, one can prove that there is 
an extremal contraction $g:\tilde{V}\to V$ such that discrepancy 
of an exceptional divisor of $g$ with respect to the pair $(V, c\H)$,
where $c<\frac{1}{\mu}$ is a canonical threshold of the pair $(V, \H)$,
is zero (such divisor is called \emph{a maximal singularity}, and 
its center on $V$ is called \emph{a maximal center}). 
One can check that there is a link of type I or II starting with 
this contraction and decreasing a properly defined ``degree'' of
the map $\chi$, all arising varieties have only terminal singularities and 
all divisorial contractions (as well as nonbirational ones) are extremal 
contractions in the sense of the ordinary Minimal Model Program
(see~\cite{Corti} or~\cite{Matsuki} for details).

The previous statements imply the following: to prove that 
$V$ cannot be transformed to another Mori fibration (i.\,e. is birationally 
rigid) it suffices to check that there are no maximal centers on $V$ except 
those that are associated with birational automorphisms of $V$, 
and to describe all birational selfmaps $\chi:V\dasharrow V$
it is sufficient to classify all maximal centers. 
In the next two sections we are going to do this for the variety $X$ described 
in section~\ref{section:setup}. 
We are also going to use the notations introduced in the current section 
throughout the paper.

\section{Excluding maximal centers: points}
\label{section:max-points}

The absence of the points that are maximal centers
is easily implied by the standard statements.

\begin{theorem}[$4n^2$-inequality, 
see~\cite{Pukhlikov-essentials} or~\cite{Corti}]
Let $V$ be a three-dimensional variety such that the linear system
$|-K_V|$ is free; let $\H\subset|-\mu K_V|$. Consider a smooth point
$x\in V$ and an intersection $Z=D_1D_2$ of two general divisors from $\H$. 
If $x$ is a maximal center (for $\H$), then $\mult_x Z > 4\mu^2$.
\end{theorem}

\begin{corollary}
\label{corollary:max-smooth-point}
A smooth point cannot be a maximal center on $X$.
\end{corollary}
\begin{proof}
Indeed, if a smooth point $x$ is a maximal center, then for general 
$D_1, D_2\in \H$ we have
$$4\mu^2<\mult_x D_1D_2\le K_XD_1D_2=\mu^2\deg X=4\mu^2,$$
a contradiction.
\end{proof}

To exclude singular points we'll need the following theorem.

\begin{theorem}[see., for example,~\cite{Cheltsov-survey}]
\label{theorem:max-sing-point}
Let $\dim V\ge 3$; let $x\in V$ be a simple double point and
$D$\t an effective ($\Q$-)divisor such that the pair $(V, D)$ is not 
canonical at $x$. Then $\mult_x D>1$.
\end{theorem}

\begin{corollary}
\label{corollary:max-sing-point}
A singular point cannot be a maximal center on $X$.
\end{corollary}
\begin{proof}
Assume that a singular point $x\in X$ is a maximal center; 
recall that all singular points of $X$ are simple double points.
Let $g:\nolinebreak\tilde{X}\to\nolinebreak X$ 
be a blow-up of the point $x$ with an exceptional 
divisor $E$; let $L\subset Q$ be a general line passing through
$\phix(x)$, $L'=\phix^{-1}(L)$, and $\tilde{L'}$\t a preimage of 
$L'$ under the map $g$. Note that $L'$ is singular at $x$, and 
$E\tilde{L'}=2$. Let $D\in\frac{1}{\mu}\H$ be a general divisor, 
$g^*(\mu D)=g_*^{-1}(\mu D)+\nu E$. By Theorem~\ref{theorem:max-sing-point} 
we have $\nu>\mu$. Hence
$$2\mu=\mu D L'=g^*(\mu D)\tilde{L'}\ge \nu E\tilde{L'}\ge 2\nu > 2\mu,$$
a contradiction.
\end{proof}

\begin{remark}
Similar statements were proved in~\cite{Grinenko}, but
the references to $4n^2$-inequality and Theorem~\ref{theorem:max-sing-point}
make the proofs much easier.
\end{remark}

\section{Excluding maximal centers:curves}
\label{section:max-curves}

We'll use the following description of low-degree curves 
on the variety $X$ (recall that the degrees are always calculated 
with respect to the anticanonical linear system).

\begin{lemma}[see, for example,~\cite{Iskovskikh-rigid}]
\label{lemma:small-degree-curves}
Let $Z\subset X$ be a curve not contained in the ramification divisor,
$\deg Z\le 3$, $\phix(Z)=Y\subset Q$.
Then one of the following cases occur.

1. $\deg Z=1$, $Y$ is a line tangent to $S$ at two (possibly coinciding)
points (it is also possible that in one or both of these points 
is just a singular point of $S$), $\left.\phix\right|_Z:Z\to Y$ 
is an isomorphism, 
$\phix^{-1}(Y)=Z\cup\delta(Z)$.

2. a) $\deg Z=2$, $Z$ is a smooth curve, $p_a(Z)=1$, $Y$ is a line intersecting
$S$ at two different points.

b) $\deg Z=2$, $Z$ is a curve with one double point, 
$p_a(Z)=1$, $Y$ is a line tangent to $S$ at a single point.

c) $\deg Z=2$, $Z$ is a smooth curve, 
$Y$ is a conic tangent to $S$ at four points 
(some of them may coincide), 
$\left.\phix\right|_Z:Z\to Y$ is an isomorphism, 
$\phix^{-1}=Z\cup\delta(Z)$.

3. $\deg Z=3$, $Y$ is a twisted cubic, $\left.\phix\right|_Z:Z\to Y$ is an 
isomorphism.
\end{lemma}

Let us describe an involution associated to a line on $X$ not contained in the 
ramification divisor.

\begin{example}[see.~\cite{Iskovskikh-rigid} or~\cite{Grinenko}]
\label{example:involution}
Let $B\subset X$ be a line not contained in the ramification divisor 
(as in case~1 of Lemma~\ref{lemma:small-degree-curves}), $B'=\delta(B)$. 
The linear system $|-K_X-B-B'|$ gives a rational map  
$\psi:X\dasharrow\P^2$; its general fiber is an elliptic curve\t
a preimage of a general line on $Q$ intersecting the line $\phix(B)$. 
This map may be regularized by passing to a variety $\tilde{X}$ that is 
a subsequent blow-up of $B$ and of a strict transform of $B'$. 
The preimage $E'$ of the curve $B'$ gives a section 
of the elliptic fibration $\tilde{\psi}:\tilde{X}\to\P^2$. Reflection
with respect to this section gives a biregular involution on an open subset 
of $\tilde{X}$ that gives rise to a biregular involution of $\tilde{X}$ 
and a birational involution of $X$. 
The action of this involution on $\Pic(X)$ is computed 
in~\cite{Iskovskikh-rigid}. Note that unlike the involutions
of a singular quartic (cf.~\cite{Mella} or~\cite{Pukhlikov-quartic}) 
the involution $\tau_B$ is constructed uniformly
regardless to the number of singular points on $B$ (moreover, there 
may be no singular points on $B$ at all).
\end{example}

From the point of view of Sarkisov program the involution $\tau_B$
is a link of type II such that any decomposition of the map 
$\chi$ starts with $\tau_B$ 
provided that the curve $B$ is a maximal center for $\H$.
Now we'll check that no other curves except the lines not contained in the
ramification divisor cannot be maximal centers. Since neither smooth
nor singular points of $X$ can be maximal centers 
(see Corollary~\ref{corollary:max-smooth-point} and 
Corollary~\ref{corollary:max-sing-point}), 
according to section~\ref{section:preliminaries} this will imply
that there are no birational maps between $X$ and other Mori 
fibrations, and that any birational selfmap $\chi:X\dasharrow X$ can 
be decomposed into a composition of the maps described in
Example~\ref{example:involution}.

\begin{lemma}
Let $Z\subset X$ be a curve not contained in the ramification divisor, and
$\deg Z>1$. Then $Z$ is not a maximal center.
\end{lemma}
\begin{proof}
Let multiplicity of $\H$ in the curve $Z$ be equal to $\nu$. Assume that 
$\nu>\mu$. By a standard computation we get $\deg Z<\deg X=4$.
Consider several possible cases depending on what the curve $Z$ is 
(we label them as in Lemma~\ref{lemma:small-degree-curves}). 
By $H$ we'll always denote a general member of the anticanonical system.

\emph{Case 2a)} 
The curve $Z$ does not pass through the singular points of $X$, and we
may use, for example, the following computation from~\cite{Iskovskikh-rigid}. 
Let $f:X'\to X$ be a blow-up of the curve $Z$, $E=f^{-1}(Z)$. Then
$$0\le (\mu f^*H-\nu E)^2(f^*H-E)=4\mu^2-4\mu\nu,$$
a contradiction.

\emph{Case 2b)} The curve $Z$ either does not pass through singular points 
of $X$, or passes through only one of them, and this point is the singular 
point $x_0$ of $Z$. Let $g:\tilde{X}\to X$ be a blow-up of the point
$x_0$, $f:X'\to\tilde{X}$\t a blow-up of a strict transform
$\tilde{Z}$ of the curve $Z$, $E_0=g^{-1}(x_0)$, $E=f^{-1}(\tilde{Z})$.
Let us define the multiplicity $\nu_0=\mult_{x_0}\H$ of the linear 
system $\H$ in the (possibly singular) point $x_0$ by an equality 
$$g^*\H=g_*^{-1}\H+(\mult_{x_0}\H)E_0;$$
if $x_0$ is a smooth point this is the multiplicity in the usual sense. 
If $x_0$ is smooth on $X$, then
\begin{multline*}
0\le(\mu (fg)^*H-\nu_0f^*E_0-\nu E)^2((fg)^*H-f^*E_0-E)=\\
=4\mu^2-4\mu\nu-(\nu_0-2\nu)^2<0,
\end{multline*}
a contradiction. If $x_0$ is a singular point of $X$, then
\begin{multline*}
0\le(\mu (fg)^*H-\nu_0f^*E_0-\nu E)^2((fg)^*H-f^*E_0-E)=\\
=4\mu^2-4\mu\nu-2(\nu_0-\nu)^2<0,
\end{multline*}
a contradiction.

\emph{Case 2c)}
Let $Z'=\delta(Z)$, and $\mult_{Z'}\H=\nu'$. Assume that 
$Z$ passes through $k$ singular points of $X$, $0\le k\le 4$ (note that these
points must be contained in $Z\cap Z'$). 
Let $U$ be a general anticanonical divisor passing through $Z$ 
(and hence through $Z'$), and $\left.\H\right|_U=\nu Z+\nu'Z'+C$.
Then on the surface $U$ the following equalities hold: 
$Z'^2=-2+\frac{k}{2}$, $Z'Z=4-\frac{k}{2}$. Hence for a general
$D\in\H$ we have
$$2\mu=DZ' \ge (4-\frac{k}{2})\nu-(2-\frac{k}{2})\nu',$$
that implies $\nu'>\nu>\mu$. The latter contradicts the fact that 
the degree of a curve that is a maximal center on $X$ cannot exceed $3$.

\emph{Case 3.}
Let $Z$ pass through $k$ singular points of $X$, say, $x_1, \ldots, x_k$.
Since $Z$ is not contained in the ramification divisor, we have $0\le k\le 6$. 
Let 
$g:\tilde{X}\to X$ be a blow-up of the points $x_1, \ldots, x_k$,
$f:X'\to\tilde{X}$\t a blow-up of a strict transform $\tilde{Z}$ of the curve $Z$, $E_i=g^{-1}(x_i)$, $E=f^{-1}(\tilde{Z})$, $\nu_i=\mult_{x_i}\H$.
Then
\begin{multline*}
0\le (\mu(fg)^*H-f^*\sum\nu_iE_i-\nu E)^2(2(fg)^*H-f^*\sum E_i-E)=\\
=8\mu^2-6\mu\nu-5\nu^2-2\sum\nu_i^2+2\nu\sum\nu_i=\\
=(8\mu^2-6\mu\nu-2\nu^2)-\frac{1}{2}(4\sum\nu_i^2-4\nu\sum\nu_i+k\nu^2)
-(3-\frac{k}{2})\nu^2=\\
=(8\mu^2-6\mu\nu-2\nu^2)-\frac{1}{2}\sum(2\nu_i-\nu)^2
-(3-\frac{k}{2})\nu^2<0,
\end{multline*}
a contradiction.
\end{proof}

\begin{lemma}[see~\cite{Iskovskikh-rigid}, \cite{CheltsovPark}]
Let $Z\subset X$ be a curve contained in the ramification divisor. 
Then $Z$ is not a maximal center.
\end{lemma}
\begin{proof}
Assume that $\mult_Z\H>\mu$. 
Let $p$ be a general point of the curve $Z$; let $L\subset Q$ be a line 
tangent to $S$ at the point $\phix(p)$, $\tilde{L}=\phix^{-1}(L)$. Then
$\tilde{L}$ is singular at $p$, and generality of $p$ implies that 
$\tilde{L}\not\subset \Bs\H$. Hence for a general $D\in\H$ we have
$$2\mu=\mu\deg\tilde{L}=\tilde{L}D\ge\mult_p\tilde{L}\mult_p\H\ge
2\mult_Z\H>2\mu,$$
a contradiction.
\end{proof}

\section{Relations}
\label{section:relations}

To prove that involutions $\tau_B$ are independent in the group $\Bir(X)$ 
it suffices to check that two lines $Z_1$ and $Z_2$ cannot appear
simultaneously as maximal centers. Let $\mult_{Z_1}\H=\nu_1$, 
$\mult_{Z_2}\H=\nu_2$. Assume that $\nu_1>\mu$, $\nu_2>\mu$; we are going
to obtain a contradiction in all the possible cases 
(cf.~\cite{Iskovskikh-rigid} and~\cite{IskovskikhPukhlikov}).

\begin{lemma}
The lines $Z_1$ and $Z_2$ cannot intersect in a single point. 
\end{lemma}
\begin{proof}
Assume that the point $x_0$ is the only intersection point of the 
lines $Z_1$ and $Z_2$. Assume in addition that there are no other 
singularities of $X$
on the curves $Z_i$ except a possible singularity at $x_0$ (other cases
are treated in a similar way but require more computations). 
Let $g:\tilde{X}\to X$ be a blow-up of the point 
$x_0$; let $f:X'\to \tilde{X}$ be a blow-up of the strict transforms
$\tilde{Z_1}, \tilde{Z_2}$ of the curves $Z_1$ and $Z_2$, $E_0=g^{-1}(x_0)$, 
$E_1=f^{-1}(\tilde{Z_1})$, $E_2=f^{-1}(\tilde{Z_2})$, $\nu_0=\mult_{x_0}\H$.
If $x_0$ is nonsingular on $X$, then 
\begin{multline*}
0\le (\mu(fg)^*H-\nu_0f^*E_0-\nu_1E_1-\nu_2E_2)^2(2(fg)^*H-f^*E_0-E_1-E_2)=\\
=8\mu^2-2\mu(\nu_1+\nu_2)+2\nu_0(\nu_1+\nu_2)-4(\nu_1^2+\nu_2^2)-\nu_0^2=\\
=(8\mu^2-2\mu(\nu_1+\nu_2)-2(\nu_1^2+\nu_2^2))-\frac{1}{2}(\nu_0-2\nu_1)^2-
\frac{1}{2}(\nu_0-2\nu_2)^2<0,
\end{multline*}
a contradiction. If $x_0$ is a singular point of $X$, then
\begin{multline*}
0\le (\mu(fg)^*H-\nu_0f^*E_0-\nu_1E_1-\nu_2E_2)^2(2(fg)^*H-f^*E_0-E_1-E_2)=\\
=8\mu^2-2\mu(\nu_1+\nu_2)+2\nu_0(\nu_1+\nu_2)-3(\nu_1^2+\nu_2^2)-2\nu_0^2=\\
=(8\mu^2-2\mu(\nu_1+\nu_2)-2(\nu_1^2+\nu_2^2))-(\nu_0-\nu_1)^2-
(\nu_0-\nu_2)^2<0,
\end{multline*}
a contradiction.
\end{proof}

\begin{lemma}
The lines $Z_1$ and $Z_2$ cannot intersect in two points (or, equivalently,
$Z_2\neq \delta(Z_1)$).
\end{lemma}
\begin{proof}
Assume that they do. Let $p$ be a general point of the curve
$Z_1$, $L\subset Q$\t a general line passing through the point $\phix(p)$, 
$\tilde{L}=\phix^{-1}(L)$. 
Then $L\cap Z_1=p$, $L\cap Z_2=\delta(p)$, and 
$\tilde{L}\not\subset \Bs\H$. Hence for a general $D\in\H$ we have
$$2\mu=\mu\deg\tilde{L}=\tilde{L}D\ge\mult_p\H+\mult_{\delta(p)}\H
\ge \nu_1+\nu_2>2\mu,$$
a contradiction.
\end{proof}

\begin{lemma}
$Z_1\cap Z_2\neq\varnothing$
\end{lemma}
\begin{proof}
Assume that $Z_1\cap Z_2=\varnothing$.
Note that there is a one-parameter family of lines $Q$ that intersect 
$\phix(Z_1)$ and $\phix(Z_2)$. Let $L$ be a general line of this family, 
$\tilde{L}=\phix^{-1}(L)$. 
Then $\tilde{L}\not\subset \Bs\H$. Hence for a general $D\in\H$ we have
$$2\mu=\mu\deg\tilde{L}=\tilde{L}D\ge \nu_1+\nu_2>2\mu,$$
a contradiction.
\end{proof}

\section{$\Q$-factoriality}\label{section:Q-factoriality}

In this section we prove Proposition~\ref{proposition:Q-factoriality}.

Let us start with an example of non-$\Q$-factorial variety $X$.

\begin{example}\label{example:non-Q-factorial}
Let a quartic $W$ be given by an equation $f_2(x)^2+f_1(x)f_3(x)=0$, 
where $f_i$ is a (general) form of degree $i$. Then the preimage of the divisor
$(f_1=0)\subset Q$ under the map $\phix$ 
splits into a union of two divisors of degree $2$ that implies 
non-factoriality of $X$ (and hence non-$\Q$-factoriality as well, since
the latter is equivalent to the former in the case of nodal
threefolds).

The variety $X$ has $12$ singular points: these are $12$ intersection 
points of the hypersurfaces $f_1=0$, $f_2=0$, $f_3=0$ and $Q$ in $\P^4$. 

Note that this variety may be also described in a different way. 
Consider a cone $K\subset\P^5$ with a vertex $P$ over a smooth 
three-dimensional quadric $Q\subset\P^4$
and a general cubic $C$ passing through $P$. Let $X'=K\cap C$; let
$\pi:X'\dasharrow Q$ be a projection of $X'$ from the point $P$. 
The birational map $\pi$ gives rise to a birational morphism 
$\tilde{\pi}=\pi\circ\sigma:\tilde{X}\to Q$, where $\sigma:\tilde{X}\to X'$ 
is a blow-up of $X'$ at the point $P$; the latter morphism is a composition 
$\tilde{X}\stackrel{\phi}{\to}X\stackrel{\varphi}{\to} Q$,
where $\phi$ is a contraction of $12$ lines on $X'$ passing through $P$
and $\varphi$ is a double cover. 
The variety $X$ has $12$ simple double singularities;
since $\phi$ is a small contraction, $X$ is not $\Q$-factorial.
\end{example}

\begin{remark}[{cf.~\cite[Example~6]{Mella}
and~\cite[Example~1.21]{Cheltsov-points}}]
Consider a general complete intersection $Y$ of a quadric $K$ and a cubic $C$ 
in $\P^5$, such that $Y$ has one simple double singularity $P$. 
Then $Y$ containes $12$ lines $l_1, \ldots, l_{12}$
passing through $P$; all these lines are contained in a common 
tangent space $\tilde{L}$ to $K$ and $C$ at $P$ (the generality 
condition lets us assume that both $K$ and $C$ are smooth). 
A projection $\pi:Y\dasharrow\P^4$ from the point $P$
is a composition of a blow-up $\sigma:\tilde{Y}\to Y$ of the point $P$ 
and a contraction
$\phi:\tilde{Y}\to Z$ of the preimages of the $12$ lines. It is
easy to see that the image $Z$ of the variety $Y$ under this projection
is a (non-$\Q$-factorial) quartic, containing a quadric surface\t
an image of an exceptional divisor of the blow-up\t
and having simple double singularities at the points $\pi(l_i)$.

One can check that the quartic $Z$ containes two quadric surfaces
and has two (projective) small resolutions;
one of them is the variety $\tilde{Y}$ described above, and the other 
is a variety $\tilde{Y'}$ that is also obtained as a blow-up 
of a singular complete intersection $Y'$ of a quadric and a cubic in $\P^5$
in its unique singular (simple double) point $P'$ (see~\cite[Example~6]{Mella}
or~\cite[Example~1.21]{Cheltsov-points}). The varieties $Y$ and $Y'$
are not isomorphic in general, but there is a natural
birational map $\tilde{\Delta}:Y\dasharrow Y'$.

When the quadric $K$ degenerates to a cone over a nonsingular three-dimensional
quadric $Q$, and the cubic $C$ passes through a vertex $P$ of the cone $K$,
the variety $Z$ obtaines a structure of a double cover of the quadric 
$Q$ (see the second description in Example~\ref{example:non-Q-factorial});
the corresponding variety $Y$ has a birational involution 
$\tilde{\delta}:Y\dasharrow Y$ arising from an involution of a double cover
$\delta:Z\to Z$.

Note also that a simpler example of a non-$\Q$-factorial quartic
with simple double singularities\t a general quartic containing a plane\t
does not degenerate to a non-$\Q$-factorial double quadric
with the same number of simple double points.
\end{remark}

Now let the divisor $S$ be singular at the points $p_1, \ldots, p_s$, $s<12$. 

In the case of a nodal threefold $\Q$-factoriality is equivalent to
factoriality; on the other hand, factoriality of a nodal Fano threefold
$X$ is equivalent to a topological condition $\rk H^2(X, \Z)=\rk H_4(X, \Z)$
(in our case $\rk H^2(X, \Z)=\rho(X)=1$).
Hence to prove $\Q$-factoriality of $X$ it suffices to check that
for a small resolution $h:\tilde{X}\to X$ we have
$\rho(\tilde{X})=s+1$. To check this it is sufficient
to show that, using the terminology of~\cite{Cynk},
the defect of $X$ equals zero (cf. the proof of Theorem~2
in~\cite{Cynk}, and also~\cite{CheltsovPark}). The latter condition
means that the singular points $p_1, \ldots, p_s$ impose independent
conditions on the hypersurfaces of degree $3$ in $\P^4$, i.\,e. that
for any point $p_i$, $1\le i\le s$, there is a cubic hypersurface
$D_i\subset\P^4$ such that $p_i\not\in D_i$ and $p_j\in D_i$ for $j\neq i$.
In the remaining part of the section we check this condition for our 
variety~$X$.

We may assume $s=11$.
We'll need the following theorem, proved in~\cite{EisenbudKoh}.

\begin{theorem}\label{theorem:Eisenbud}
The points $p_1, \ldots, p_s\in\P^n$ impose independent conditions on forms
of degree $d$ if each linear subspace of dimension $k$ contains
at most $dk+1$ of the points $p_1, \ldots, p_s$.
\end{theorem}

In the remaining part of the section we'll check the assumptions of
Theorem~\ref{theorem:Eisenbud} and derive
Proposition~\ref{proposition:Q-factoriality} from 
Theorem~\ref{theorem:Eisenbud}.

\begin{lemma}[{see~\cite[Corollary~2.5]{Shramov}}]
\label{corollary:10-on-a-quadric-in-p3}
Let $Y\subset\P^3$ be an irreducible quadric; let
$p_1, \ldots, p_{10}, q\in Y$ be such points that
no line contains $4$ of the points $p_1, \ldots, p_{10}$,
no conic contains $7$ and no twisted cubic containes $10$ of them.
Then there is a divisor $D\in\left.\O_{\P^3}(3)\right|_{Y}$ passing 
through $p_1, \ldots, p_{10}$ and not passing through $q$.
\end{lemma}

\begin{lemma}\label{lemma:4-on-a-line}
A line contains at most $3$ points of $p_1, \ldots, p_{11}$. 
\end{lemma}
\begin{proof}
Assume that $4$ points, say, $p_1, \ldots, p_4$, lie on a line $l$. 
Then $l\subset Q$, and since $W$ is tangent to $Q$ at $p_1, \ldots, p_4$,
it follows that $W$ is tangent to $l$ at $p_1, \ldots, p_4$, and hence
$l\subset W$. Choose homogeneous coordinates  
$(x_0:\ldots: x_4)$ so that $Q$ is given by an equation 
$1/2(x_0^2+\ldots+x_4^2)=0$, and $l$ is given by equations $x_2=x_3=x_4=0$. 
Then $W$ is defined by an equation  
$$x_2C_2(x)+x_3C_3(x)+x_4C_4(x)=0,$$
where $\deg C_i=3$. The divisor $S$ is singular at the points of $l$ where
the matrix of partial derivatives 
$$
\left(
\begin{array}{ccccc}
x_0 & x_1 & x_2 & x_3 & x_4 \\
0 & 0 & C_2 & C_3 & C_4 
\end{array}
\right)
$$
has rank $1$. In these points the minors 
$x_0C_2$ and $x_1C_2$ must vanish. Since there are at least 
$4$ such points, $C_2(x)$ is identically zero along $l$.
Similar arguments show that $C_3$ and $C_4$ also vanish along $l$, 
and $S$ has non-isolated singularities, that contradicts our assumptions.
\end{proof}

\begin{lemma}\label{lemma:7-on-a-plane}
A plane contains at most $6$ of the points $p_1, \ldots, p_{11}$.
\end{lemma}
\begin{proof}
Assume that $7$ points, say, $p_1, \ldots, p_7$, lie in a plane $L$. 
If $L$ intersects $Q$ along a reducible conic, 
then there are $4$ collinear points among $p_1, \ldots, 
p_7$ that contradicts Lemma~\ref{lemma:4-on-a-line}. Hence the conic
$L\cap Q$ is nonsingular, and we may assume that $Q$ is given by an equation
$1/2(x_0^2+\ldots+x_4^2)=0$, $L$ is given by equations $x_3=x_4=0$, 
and $W$\t by an equation
$$x_3C_3(x)+x_4C_4(x)+Q(x)\tilde{Q}(x)=0,$$
where $\deg C_i=3$, $\deg\tilde{Q}=2$, and  $Q(x)=0$ is an equation of
the quadric $Q$.
The divisor $S$ is singular at the points of $L\cap Q$ where the matrix
of partial derivatives  
$$
\left(
\begin{array}{ccccc}
x_0 & x_1 & x_2 & x_3 & x_4 \\
x_0\tilde{Q} & x_1\tilde{Q} & x_2\tilde{Q} & C_3+x_3\tilde{Q} & 
C_4+x_4\tilde{Q} 
\end{array}
\right)
$$
has rank $1$. In these points the minors
$x_0(C_3+x_3\tilde{Q})-x_3x_0\tilde{Q}=x_0C_3$ and 
$x_1(C_3+x_3\tilde{Q})-x_3x_1\tilde{Q}=x_1C_3$ must vanish. Since there are
at least $7$ such points, $C_3$ and, similarly, $C_4$ must vanish along
$L\cap Q$; hence $S$ has non-isolated singularities, that is a contradiction.
\end{proof}

\begin{lemma}
\label{lemma:10-on-a-twisted-cubic}
A twisted cubic containes at most $9$ of the points $p_1, \ldots, p_{11}$.
\end{lemma}
\begin{proof}
Assume that there is a twisted cubic $\Gamma$ containing $10$ of the points 
$p_1, \ldots, p_{11}$. The twisted cubics contained in $Q$ form two orbits 
with respect to the action of $\Aut(Q)=\PSO(5)$: one of them consists of
the twisted cubics $\gamma$ such that their linear span $L(\gamma)$ intersects 
$Q$ by a nonsingular quadric surface, and the other consists of the 
twisted cubics $\gamma$ such that $L(\gamma)$ intersects $Q$ by a quadric cone.
In particular, we may assume that in appropriate homogeneous coordinates 
$(x_0:\ldots:x_4)$ the curve $\Gamma$
is given by equations 
\begin{gather*}
x_1^2=x_0x_2, \quad x_2^2=x_1x_3,\\
x_0x_3=x_1x_2, \quad x_4=0,
\end{gather*}
and $Q$ is given either by equation $x_0x_3-x_1x_2+x_4^2=0$,
or by equation $x_1^2-x_0x_2+x_3x_4=0$. The quartic $W$ is given by 
$$(x_0x_3-x_1x_2)Q_1+(x_1^2-x_0x_2)Q_2+(x_2^2-x_1x_3)Q_3+x_4C=0,$$
where $\deg Q_i=2$, $\deg C=3$. The divisor $S$ is singular 
at the points of $\Gamma$ where the matrix $M$ of partial derivatives 
has rank $1$. In either case one can check that if $M$ has rank $1$ at $10$ 
points on $\Gamma$, then it has rank $1$ along $\Gamma$, i.\,e. $S$ has 
non-isolated singularities, a contradiction.
\end{proof}

\begin{lemma}\label{lemma:10-in-a-space}
Assume that $10$ of the points $p_1, \ldots, p_{11}$ are contained 
in a three-dimensional space. 
Then the points $p_1, \ldots, p_{11}$ impose independent conditions
on the hypersurfaces of degree~$3$.
\end{lemma}
\begin{proof}
Let $p_1, \ldots, p_{10}\in H\simeq\P^3$. If $p_{11}\in H$, 
we apply Lemma~\ref{corollary:10-on-a-quadric-in-p3} to a quadric
$Y=Q\cap H$ (choosing different points $p_i$ as the distinguished point $q$; 
the assumptions of Lemma~\ref{corollary:10-on-a-quadric-in-p3} hold
due to Lemma~\ref{lemma:4-on-a-line}, Lemma~\ref{lemma:7-on-a-plane}
and Lemma~\ref{lemma:10-on-a-twisted-cubic}).
If $p_{11}\notin H$, then we may choose a cone $K$ over a divisor 
$D_{11}\in\left.\O(3)\right|_Y$ passing through $p_1, \ldots, p_{10}$,
so that $p_{11}\not\in K$. On the other hand, 
Lemma~\ref{corollary:10-on-a-quadric-in-p3} implies that for any point
$p_i$, $1\le i\le 10$, there is a divisor $D_i\in\left.\O(3)\right|_Y$
passing through all the points $p_1, \ldots, p_{10}$ except $p_i$; a cone 
with a vertex $p_{11}$ over such divisor passes through all the points 
$p_1, \ldots, p_{11}$ except $p_i$.
\end{proof}

\begin{proof}[Completion of the proof of~\ref{proposition:Q-factoriality}]
Let us check that the points $p_1, \ldots, p_{11}$ impose independent 
conditions on the forms of degree $3$ in $\P^4$.
In the notations of Theorem~\ref{theorem:Eisenbud} we have $n=4$, $d=3$. 
No $4$ of the points $p_1, \ldots, p_{11}$ are collinear by 
Lemma~\ref{lemma:4-on-a-line}; no $7$ of them are coplanar due to
Lemma~\ref{lemma:7-on-a-plane}. By Lemma~\ref{lemma:10-in-a-space} we may  
assume that no $10$ of the points $p_1, \ldots, p_{11}$ lie in 
a three-dimensional space. Hence the assumptions of 
Theorem~\ref{theorem:Eisenbud} hold, and we are done.
\end{proof}

\end{document}